\documentclass{article}

\usepackage{amssymb,amsmath,amsthm}

\newtheorem*{Theorem}{Theorem}
\newtheorem*{Cor}{Corollary}
\newtheorem{proposition}{Proposition}
\newtheorem*{lemma}{Lemma}
\newtheorem*{claim}{Claim}

\begin{document}

\title
% [on minimality of foliations]
{A note on minimality of foliations \\
for partially hyperbolic diffeomorphisms}

\author{Katsutoshi SHINOHARA}

%% \address{3-8-1 Komaba, Meguro, Tokyo, 153-8914 Japan}
%% \email{sinon@ms.u-tokyo.ac.jp}
%% \subjclass[2000]
%% {Primary 37D30}
%% \keywords{partial hyperbolicity, robust transitivety}

\maketitle

\begin{abstract}
It was shown that in robustly transitive, 
partially hyperbolic diffeomorphisms on three dimensional closed  manifolds,
the strong stable or unstable foliation is minimal.
In this article, we prove ``almost all'' leaves of both 
stable and unstable foliations are dense in the whole manifold.
\end{abstract}

\section{Introduction}
In \cite{BDU}, the minimality of strong stable 
and unstable foliations in robustly transitive partially hyperbolic 
diffeomorphisms is studied.
When investigating dynamics from the measure theoretic point of view, 
the minimality of stable and unstable foliations often plays a crucial role,
so the study is an interesting problem.

One of the results in \cite{BDU} is the following: 
in the robustly transitive, partially hyperbolic diffeomorphisms 
on three dimensional closed  manifolds,
the strong stable  {\it or} unstable foliation is minimal.
Note that the minimality of {\it both} stable and unstable foliation
is still an open problem.
Roughly speaking, this asymmetry comes from 
the nature of center foliation.

In this article, applying the method in \cite{BDU},
we prove ``almost all'' leaves of both 
stable and unstable foliations are dense in the whole manifold.

\section{Notations and Definitions}
Let $M$ denote a smooth, compact, three-dimensional 
Riemannian manifold,
${\rm Diff} ^1 ( M )$ be the space of $C^1$ 
diffeomorphisms of $M$ with $C^1$ topology.

A diffeomorphism $f\colon M\to M$ is called {\it transitive}
if there exists $x \in M $ such that 
$ \{ f^n (x) \mid n \ge 0 \} $ is dense in $M$.
Note that this condition is equivalent to the following: if $U$ and $V$
are nonempty open sets in $M$, then there exists $n \ge 0$
such that $U \cap f^n(V) \neq \emptyset$. 
With this equivalence, one can deduce 
that $f$ is transitive if and only if $f^{-1}$ is .

A diffeomorphism $f$ is called 
{\it robustly transitive} if there exists a neighborhood $\mathcal{U}$ 
($\subset {\rm Diff} ^1 ( M ) $) of $f$ 
such that for all $g \in \mathcal{U}$, $g$ is transitive.
By the continuity of the map
$(\, \, \cdot \, \, ) ^{-1} \colon {\rm Diff} ^1 ( M ) \to {\rm Diff} ^1 ( M )$,
$f$ is robustly transitive if and only if $f^{-1}$ is.

Let $f$ be a diffeomorphism of $M$ and 
$\Lambda$ be an invariant subset of $M$. 
$f$ is called 
{\it strongly partially hyperbolic on $\Lambda$} if there exists
a continuous invariant splitting of tangent bundle 
$TM|_{\Lambda} = E^{s} \oplus E^{c} \oplus E^{u} $
and $0 < \lambda < 1$ which satisfy following conditions:
\begin{enumerate}
\item for all $x \in \Lambda$,
\[
\| df | _{E^{s}} (x) \| < \lambda \quad \mbox{and} \quad
\| df^{-1} | _{E^{u}} (x) \|  < \lambda ;
\]

\item for all $x \in \Lambda$, 
\[
\| df | _{E^{s}} (x) \| \| df^{-1} | _{E^{c}} (x) \| < \lambda
\quad \mbox{and} \quad
\| df | _{E^{c}} (x) \| \| df^{-1} | _{E^{u}} (x) \| < \lambda.
\]
\end{enumerate}
In the sentences above, $\| \, \cdot \, \|$ denotes the operator norm.

We focus on 
robustly transitive, strongly partially hyperbolic (on $M$), 
but not Anosov diffeomorphisms.
One can perturb such diffeomorphisms having hyperbolic periodic points
with index  1 and 2 
(where the index of hyperbolic periodic point is the dimension of its stable manifold).
This is a consequence of the ergodic closing lemma of Ma\~{n}\'{e}, see \cite{Ma2}.

We denote by $\mathcal{RP}(M)$ the set of diffeomorphisms on $M$
which is robustly transitive,
strongly partially hyperbolic,
and has hyperbolic periodic points with index 1 and 2.
In ${\rm Diff} ^1 ( M )$, the set $\mathcal{RP}(M)$ is open.
The example of Ma\~{n}\'{e} (see \cite{Ma1}) shows $\mathcal{RP}(T^3) \neq \emptyset$.

\section{Statement of the main theorem}
Let $f\in \mathcal{RP}(M)$. The uniform hyperbolicity of 
$E^{s}$ (resp. $E^{u}$) gives the
{\it strongly stable foliation}  $\{ \mathcal{F}^{s}(x) \}$ 
(resp. {\it strongly unstable foliation} $\{ \mathcal{F}^{u}(x) \}$). 
See \cite{HPS} for detail.

A foliation is called {\it minimal} if every leaf is dense in $M$.
We are interested in whether $\{ \mathcal{F}^{s}(x) \}$ and 
$\{ \mathcal{F}^{u}(x) \}$ are minimal or not.
The following is the main theorem in this paper.

\begin{Theorem}
There exists an open dense subset 
$\mathcal{U}^{s}_{0}$ (resp.\,$\mathcal{U}^{u}_{0}$) of $\mathcal{RP}(M)$,
such that for all $g \in \mathcal{U}^{s}_{0}$ (resp.\,$\mathcal{U}^{u}_{0}$) 
there exists a residual subset
 $\mathcal{O}^{s} (g)$ (resp. $\mathcal{O}^{u} (g)$) of $M$, 
satisfying the following property: 
if a point $p$ is in $\mathcal{O}^{s} (g)$ (resp. $\mathcal{O}^{u} (g)$) 
then the leaf $\mathcal{F}^{s}(p)$ (resp. $\mathcal{F}^{u}(p)$)
is dense in $M$.
\end{Theorem}

Then we have the following corollary:

\begin{Cor}
There exists an open dense subset 
$\mathcal{U}_{0}$ of $\mathcal{RP}(M)$,
such that for all $g \in \mathcal{U}_{0}$ 
there exists a residual subset
$\mathcal{O}(g)$ satisfying the following property: 
if $p \in \mathcal{O}(g)$ 
then the leaf $\mathcal{F}^{s}(p)$ and $\mathcal{F}^{u}(p)$
are dense in $M$.
\end{Cor}

One can get $\mathcal{U}_{0}$ as the intersection of 
$\mathcal{U}_{0}^{s}$ and $\mathcal{U}_0^{u}$,
and for $g \in \mathcal{U}_{0}$ the open set 
$\mathcal{O}(g)$ is attained as the intersection of 
$\mathcal{O}^{s} (g)$ and $\mathcal{O}^{u} (g)$.

In the next section, we prove only the statement about the stable foliation.
The proof of the unstable one can be given by applying our argument to $f^{-1}$.
  
\section{Proof of the theorem}
First let us recall the following two propositions 
in \cite{BDU}:

\begin{proposition}[\cite{BDU}, Theorem 2.1]
\label{BDU1}
Let $f \in \mathcal{RP}(M)$ and $p$ (resp. $q$) be 
a hyperbolic periodic point of index $1$ (resp. $2$).
Then the unstable manifold $W^{u}(p)$(resp. stable manifold $W^{s}(q)$) is dense in $M$.
\end{proposition}

\begin{proposition}[\cite{BDU}, Theorem 2.6]
\label{BDU2}
Let $f$ and $p$ be as above and $\mathcal{V}_{f}$ be the open neighborhood
of f in which the continuation $p_g$ of $p$ can be defined.
Then, there exists an open dense subset $\mathcal{W}_{f}$ of $\mathcal{V}_{f}$
such that for all $g \in \mathcal{W}_{f}$ the stable manifold of 
$W^{s}(p_{g})$ is dense in M.
\end{proposition}

The proof of proposition \ref{BDU1} is easy, and that of 
proposition \ref{BDU2} needs the connecting lemma. 
%First, by the density of $W^{u}(p)$ and $W^{s}(q)$, 
%one can surmise that they have a nonempty transversal intersection. 
%Secondly, by the existence of transitive orbit, one can 
%perturb $f$ so that $W^{s}(p)$ and $W^{u}(q)$ have 
%quasi-transversal intersection  perturbation, and make it arbitrary small.
%Finally, under the existence of {\it heterodimensional cycle}, 
%one can perturb $f$ so that 
%$W^{u}(p) \subset \overline{W^{u}(q)}$ and
%$W^{s}(q) \subset \overline{W^{s}(p)}$ in a robust way, 
%see \cite{DR} (note that this dynamics is far from tangency, 
%for the existence of strong partial hyperbolicity). 
%By the density of $W^{s}(q)$, the proposition follows.
Using these propositions, we will prove our theorem.
First, define $U^{s}_0$ as the set of $f \in \mathcal{RP}(M)$
satisfying following properties:
\begin{itemize}
\item $f$ has a hyperbolic periodic point $p$ 
with index 1,
\item the stable manifold of $p$ is dense in $M$.
\end{itemize}

By proposition \ref{BDU1} and \ref{BDU2},
$U^{s}_0$ is an open dense subset of $\mathcal{RP}(M)$. 
We focus on fixed $f \in U^{s}_0$, 
$p$ and $W^{s}_{\mathrm{loc}}(p)$. 

Next, let $A_0$ be an open ball centered at $p$ with radius $r$,
where $r$ is a sufficiently small positive number satisfying:
\[
x \in A_0 \Longrightarrow 
\mathcal{F}^{s}(x) \cap W^{u}_{\mathrm{loc}}(p) \neq \emptyset.
\]
We can take such $r$  by the $C^1$ continuity of the 
unstable foliation.
Define $A_n$ as the open ball centered at $p$ 
with radius $\displaystyle \frac{r}{n}$. 
These balls also satisfy the property above with $A_0$ 
replaced by $A_n$.

Now, put 
\[
 C_n = \bigcup ^{-\infty}_{k=0} f^{k}(A_n).
\]
This is an open subset of $M$. 
\begin{lemma}
For all $l \ge 0, \, C_l$ is an open dense subset of $M$.
\end{lemma}
\begin{proof}
% Density is a consequence of the transitivety of $f$. 
Let $x\in M$ and $O$ be any neighborhood of $x$.
By the transitivity of $f$, there exists $m \le 0$ such that 
$O \cap f^m (A_l) \neq  \emptyset$,
which means $O \cap C_l \neq \emptyset$.
Since $O$ can be arbitrary, $x \in \overline{C_l}$,
which implies $M = \overline{C_l}$.
\end{proof}
Put $\mathcal{O}^{s}(f) = \cap _{k=0}^{\infty}C_k$. 
By Baire's category theorem,
$\mathcal{O}^{s}(f)$ is a residual subset of $M$.
The next claim is nothing but our theorem.
\begin{claim}
For every point $x \in \mathcal{O}^{s}(f)$, 
$\mathcal{F}^{s} (x)$ is a dense subset of $M$.
\end{claim}
\begin{proof}
Observe that if the leaf $\mathcal{F}^{s} (x)$ passes arbitrary close to $p$,
then $\mathcal{F}^{s} (x)$ is dense in $M$
(this follows from the density of $W^{s}(p) = \mathcal{F}^{s}(p)$ 
and the continuity of $\{ \mathcal{F}^{s}(y) \}$).

Let $x \in \mathcal{O}^{s}(f)$. 
We claim that, for all $l \ge 0$, 
\[
\mathcal{F}^{s}(x) \cap (W^{u}_{\mathrm{loc} } (p) \cap A_l) \neq \emptyset.
\]

In fact,
\begin{align*}
x \in \mathcal{O}^{s}(f) 
\iff & \forall l \ge 0, \, x \in C_l \\
\iff & \forall l \ge 0, \, \exists k \le 0, \, x \in f^k(A_l) \\
\iff & \forall l \ge 0, \, \exists k \le 0, \, f^{-k} (x) \in A_l \\
\iff & \forall l \ge 0, \, \exists k \le 0, \,
\mathcal{F}^{s} \big( f^{-k}(x) \big) \cap \big(W^{u}_{\mathrm{loc}} (p) \cap A_l\big) \neq \emptyset   \\
\Longrightarrow \, \, &  \forall l \ge 0, \,  
\mathcal{F}^{s}(x) \cap \big( W^{u}_{\mathrm{loc}} (p) \cap A_l \big) \neq \emptyset.   
\end{align*}
Thus, we know if $x \in \mathcal{O}^{s}(f)$ then 
$\mathcal{F}^{s}(x)$ passes arbitrary close to $p$.
This proves the claim, and the proof of the theorem is completed.
\end{proof}


\begin{thebibliography}{9}

\bibitem{BDU} 
C. Bonatti, L. D\'{i}az, R. Ures, 
{\it Minimality of strong stable and unstable foliations for partially hyperbolic diffeomorphisms},
J. Inst. Math. Jussieu 1 (2002), no. 4, 513--541. 
\bibitem{DR}
L. D\'{i}az, J. Rocha, 
{\it Partially hyperbolic and transitive dynamics generated by heteroclinic cycles},
 Ergodic Theory Dynam. Systems 21 (2001), no. 1, 25--76.
\bibitem{HPS}
M. Hirsch, C. Pugh, M. Shub,  
{\it Invariant manifolds},  
Lecture Notes in Mathematics, Vol. 583. 
Springer-Verlag, Berlin-New York, 1977.
\bibitem{Ma1}
R. Ma\~{n}\'{e}, 
{\it Contributions to the stability conjecture},
Topology 17 (1978), no. 4, 383--396. 
\bibitem{Ma2}
R. Ma\~{n}\'{e},
{\it An ergodic closing lemma},
Ann. of Math. (2) 116 (1982), no. 3, 503--540.
\end{thebibliography}
\end{document}